\documentclass[conference]{IEEEtran}
\ifCLASSINFOpdf
\else
\fi
\usepackage{amsfonts,amssymb, amsmath, mathrsfs}
\usepackage{color}                           
\newtheorem{definition}{Definition}

\newtheorem{theorem}{Theorem}

\def\R{\mathop{\mathbb R}\nolimits}

\begin{document}
%
\title{Reduction of the Pareto Set in Multicriteria Economic Problem with CES Functions}

\author{\IEEEauthorblockN{Alexey Zakharov}
\IEEEauthorblockA{Faculty of Applied Mathematics and Control Processes\\
Saint Petersburg State University\\
7/9 Universitetskaya nab., St. Petersburg, 199034 Russia.\\
Email: a.zakharov@spbu.ru}
}

\maketitle

\begin{abstract}
A multicriteria economic problem is considered: the basic production assets
and the labor resources define a set of feasible solutions (alternatives); labor costs,
costs of the basic production assets (to be minimized),
and cost of the manufactured products (to be maximized) are objective functions.
The production function with constant elasticity of substitution is used.
The decision maker's (DM's) preferences are introduced as follows: lower labor costs and
costs of the basic production assets have the greater importance than higher income, and vice a versa.
The fuzzy preferences along with the compromise have a degree of its confidence.
Such crisp and fuzzy information is applied in the axiomatic approach of the Pareto set reduction by V.~D.~Noghin.
We show how to construct a set, which is an upper bound of the optimal choice and belongs
to the Pareto set of the problem in crisp and fuzzy cases. In fuzzy case one should solve a three crisp
multicriteria problems.
Thus, upon the crisp and fuzzy DM's preferences a narrower upper bounds of the optimal set of resources
with respect to criteria and additional information are obtained.
\end{abstract}

%
\IEEEpeerreviewmaketitle

\section{Introduction}
In economic systems the relationship between resources consumption and production output is modeled by production function.
Nowadays a class of neoclassical production functions is very widespread and
is commonly used in many researches. 
This class is defined as a class of functions satisfy the following conditions:
continuity; the output is zero, if resources are absent; the function is twice differentiable;
the first derivatives are continuous and non-negative; marginal productivity is a continuous nonincreasing function.

First who proposed a certain form of production function was Cobb and Douglas in their pioneer work~\cite{Cobb_Douglas}.
Authors empirically assigned the relationship between output volume and production resources (labor and capital)
and approximated it by the function, which now called the Cobb–-Douglas production function. It serves as a very useful
mathematical instrument in many economic research. However this function imposes {\color{black}some}
constraints on production process:
elasticity of substitution is always equaled 1. 
Another well-known production function is the Leontief production function.
Here elasticity of substitution equals 0, which means that production resources are non-substitutional.

To resolve the aforementioned restriction Arrow, Chenery, Minhas, and Solow proposed~\cite{Arrow_Chenery} a production function
with constant elasticity of substitution (CES function), where elasticity is a parameter of function and is time-invariant.
Thus, the Cobb–-Douglas and the Leontief functions are special cases of CES function.
Enumerate another models made a contribution to the theory: Mukerji function~\cite{Mukerji};
non-homotetic CES function by Sato~\cite{Sato}; translogarithmic function by Christensen, Jorgensen, Lau;
variable elasticity of substitution function (VES) by Revankar~\cite{Revankar};
generalized Leontief production function by Diewert~\cite{Diewert}.

Economic optimization problems arise as problems of profit maximization taking into account various factors and
restrictions, and researches couple with it using plenty of approaches and tools~\cite{Rasmussen}.
Nevertheless many investigated models have single criterion, but actually its complexity makes it multiple.

In our paper we assume that technology has sufficient stability with respect to the ratio of the production factors (resources).
So, in the model of economic system we take the CES function. Then we state following multicriteria choice problem. Consider two type of resources: labor and basic production assets, which compose
a set of feasible solutions (alternatives). Using the form of production function, prices of resources and output
we introduce a vector criterion with following components: the labor costs, the costs of the basic production assets,
and the cost of the manufactured products. Also we consider a binary preference relation of the DM in two cases,
crisp and fuzzy.

When solving the obtained multicriteria problem we get the Pareto set, which in practice is rather big,
and it brings difficulties in making the final optimal decision.
For that reason we introduce an additional information about the DM's preference relation in terms of so-called
''quanta of information'' (for crisp and fuzzy cases), reflecting the compromise between criteria or groups of criteria,
that the DM is ready to go.
When the preference relation is crisp using this information we establish such set of resources that composes an upper bound of the optimal choice and belongs to the initial Pareto set.

In fuzzy case the application of information comes to solving three crisp multicriteria choice problems.
The Pareto sets of these problems form a collection of nested sets of
resources assigned a fuzzy set, which is an upper bound of the ''fuzzy'' optimal choice and belongs to the initial Pareto set. Thus, we show how take into account additional information about the DM's crisp and fuzzy preference relations
in order to reduce the Pareto set. As a result, it simplifies the final optimal choice among the set of feasible resources with
respect to maximizing the cost of the manufactured products and minimizing costs of the resources.
\section{Mathematical economic model}
\label{ref_sect_crisp_model}
Traditionally optimization problem considered in economic research is a single-criterion problem
and consists in the profit maximization, as the difference between income and costs.
We are going to state the problem as a multicriteria problem taking into account income and costs separately
and moreover add to the model the preferences of the DM.
Let us start with mathematical statement of the problem.

\subsection{Mathematical background}
\label{ref_subsect_crisp_model_math}
By the formal description~\cite{Podinovskiy_Noghin} a multicriteria problem $<~X, f>$ is described by following components:
\begin{itemize}
\item a set of feasible solutions (alternatives) $X$;
\item a vector criterion $f = (f_1, f_2, \ldots,  f_m)$ defined on set $X$.
\end{itemize}

Components of the vector criterion $f_1, f_2, \ldots, f_m$ compose a goals of the DM. Under optimal solution to problem $<X, f>$
is usually interpreted the set of such solutions that can not be improved with respect to the all criteria simultaneously~\cite{Ehrgott}.
And we come to the definition of the Pareto set,
or the set of pareto-optimal solutions $P_f(X)$, which also can be interpreted as the set of non-dominated solutions
with respect to relation $\geq$~\cite{Ehrgott, Podinovskiy_Noghin}:
$$P_f(X)= \{x \in X \mid \nexists x^\ast \in X: f(x^\ast) \geq f(x) \}.$$
If we denote $Y = f(x)$, then we obtain the Pareto set $P(Y)$ with respect to the set $Y$:
$P(Y)=\{y \in Y \mid \nexists y^\ast \in Y: y^\ast~\geq~y~\}.$
Note, that the relation $y^\ast \geq y$ means that $y^\ast_i \geqslant y_i$
for all $i \in \{ 1, 2, \ldots, m \}$, and $y^\ast \neq y$.

Noghin V.D. added one more component to the mulitcriteria problem $<X, f>$
and introduced in the book~\cite{Nogin2016} a multicriteria choice problem $<X, f, \succ>$ with
\begin{itemize}
\item a set of feasible solutions (alternatives) $X$;
\item a vector criterion $f = (f_1, f_2, \ldots,  f_m)$ defined on set $X$;
\item an asymmetric binary preference relation of the DM $\succ$ defined on set $Y = f(X)$,
\end{itemize}
where the binary relation $\succ$ satisfies some axioms of so-called ''reasonable'' \ choice.
This relation expresses the wishes of the DM, and notation $f(x') \succ f(x'')$
for $x', \, x'' \in X$ means that the DM prefers the solution $x'$ to $x''$.

According to these axioms the relation $\succ$ is irreflexive, transitive,
invariant under linear positive transformation and compatible with each criteria $f_1, f_2, \ldots, f_m$.
The statement of compatibility means that the DM is interested in
increasing value of each criterion when values of others criteria are constant.
Also, we assume that when considering two possible variants $x', \ x'' \in X$
the excluded solution could not be chosen from the whole set $X$ as well.
For example, if the relation $f(x') \succ f(x'')$ holds, then solution $x''$
does not belong to the optimal choice within the whole set $X$.
Note, indeed the relation $\succ$ is defined not only on set $Y$, but on the entire criterion space~$\R^m$.

In~\cite{Nogin2016} the author establishes the {\it Edgeworth--Pareto principle}, which says that applying the axioms of ''reasonable'' choice any set of selected solutions $C(X)$ belongs to the Pareto set $P_f(X)$. Here, the set of selected solutions is interpreted as some abstract set, that satisfies all hypothetic preferences of the DM. So, the optimal choice should be done within the Pareto set only if corresponding preference relation $\succ$ fulfills the axioms of ''reasonable'' choice.
Thus, all solutions from the Pareto set are equivalent to each other with respect to the optimal choice.

In real-life multicriteria problems the Pareto set is rather big, and, obviously, the DM wants to narrow it.
So, V.~D.~Noghin proposed to use specific information on the DM's preference relation $\succ$ to reduce
the Pareto set staying within the optimal choice set~\cite{Nogin2014, Nogin2016_conf_OR, Nogin2016}.
Consider the definition of such information.

\begin{definition}
\label{ref_def_AB}
Say that there exists a {\it "quantum of information"} about the DM's preference
relation $\succ$ if the vector $y' \in \R^m$ such that
\begin{equation}
\begin{split}
\label{ref_AB}
y'_i = w_i > 0 \quad \forall i \in A, &\quad
y'_j = -w_j < 0 \quad \forall j \in B, \\
y'_s = 0 \quad \forall & s \in I\setminus (A \cup B)
\end{split}
\end{equation}
satisfies the expression $y' \succ 0_m$. Here $I = \{1, 2, \ldots, m \}$, $A, B \subset I$,
$A \neq \emptyset$, $B \neq \emptyset$, $A \cap B = \emptyset$.
In such case we will say, that the group of criteria $A$ is more important than the group of criteria $B$ with given positive parameters $w_i \ \forall i \in A$, $w_j \ \forall j \in B$.
\end{definition}

Thus, "quantum of information" shows that the DM is ready to compromise
by loosing the criterion $f_j$ by $w_j$ amount (from less important group $B$)
for gaining additional $w_i$ amount with respect to the criterion $f_i$ (from more important group $A$).

In~\cite{Nogin2016} author shows how to take into account ''quanta of information''. It consists in constructing a ''new'' vector criterion $g$ using components of ''old'' one $f$ and parameters of information $w_i$, $w_j$. Then one should find the Pareto set of ''new'' multicriteria problem with the same set of feasible solutions $X$ and a ''new'' vector criterion $g$. Obtained set $P_g(X)$ will be the reduction of the Pareto set, i.e. it will belong to the Pareto set $P_f(X)$ of the initial problem and compose a narrower upper bound of the optimal choice.

\subsection{Construction of the model}
\label{ref_econ_sys_crysp}
Now pass to the main problem of the paper. Consider an economic system of goods production using basic production assets and labor resources. Assume that someone seeks to increase its production output, reducing simultaneously the resource costs. Clearly, a higher output cannot be achieved without consuming additional resources, and so output maximization contradicts resource reduction. Such situation brings us to the optimization problem which takes into account more then one criterion.

Denote by $K$ and $L$ quantities of the basic production assets and the labor resources respectively.
According to~\cite{Arrow_Chenery} let us suppose that the relationship between the production output and
the resources consumption is modeled by the production function with constant elasticity of substitution
\begin{equation}
\label{ref_func_CES}
Q = F(a K^{-r} + (1-a)L^{-r})^{-(1/r)},
\end{equation}
where $Q$ is a quantity of output, $F$ is a factor productivity, $a$ is a share parameter, $r$ is a parameter such that
the quantity $1/(1+r)$ is the elasticity of substitution. Function~(\ref{ref_func_CES}) should satisfy the inequalities
\begin{equation}
\label{ref_func_CES_cond}
F>0, \ 0 < a <1, \ r > -1.
\end{equation}
in order to be neoclassical.

Consider the optimization problem of finding the optimal set of resources with respect to the labor costs, the costs of the basic production assets and the cost of the manufactured products. Obviously, the first and the second criteria should be minimized, the third one should be maximized. So we have more than one objective function, that proceed us to a multicriteria problem.

Using the notation of previous subsection let us build a mathematical model of aforementioned economic problem.

The set of feasible solutions $X$ consists of all pairs $(K, L)$ such that $K, L > 0$.
The objective functions (criteria) $f~=~(f_1, f_2, f_3)$ are the labor costs $f_1$,
the costs of the basic production assets $f_2$, and the cost of the manufactured products $f_3$.
Based on the aforesaid, functions $f_1$ and $f_2$ should be minimized,
and so we take them with minus sign to match the compatibility of the relation $\succ$.
Consequently, the criteria are defined by
\begin{equation}
\begin{split}
\label{ref_f}
f_1(K, L) = -p_{_K} K,& \quad f_2(K, L) = -p_{_L} L, \\
f_3(K, L) = p_{_Q} Q = p_{_Q} &F(a K^{-r} + (1-a)L^{-r})^{(-1/r)},
\end{split}
\end{equation}
where the quantities $p_{_K}$, $p_{_L}$, and $p_{_Q}$ are the prices of the corresponding resources and manufactured products.

The set $Y = f(X)$ describes a concave surface in $\R^3$, because its implicit representation is the following:
\begin{multline*}
Y = \left\{ y \in \R^3 \mid y_3 = p_Q F\left(a \left(-\frac{y_1}{p_K}\right)^{-r} + \right. \right. \\
+ \left. \left. (1-a)\left(-\frac{y_2}{p_L}\right)^{-r}\right)^{(-1/r)} \right\}.
\end{multline*}

Following the model of multicriteria choice problem $<~X, f, \succ>$ introduce an irreflexive, transitive,
and invariant under linear positive transformation binary relation $\succ$ reflecting preferences of the DM.
The condition of compatibility with respect to all components $f_1$, $f_2$, and $f_3$ is valid
due to negative sign of resources costs.
Also we suppose that from two possible variants $x', \ x'' \in X$
the excluded solution could not be chosen from the whole set $X$ as well.

Thus the Edgeworth--Pareto principle is valid for this problem,
and we should search optimal solutions within the Pareto set, and only within it.

Due to the concavity of set $Y$ we get that the Pareto set $P_f(X) = X$.
It means that any pair of resources $(K, L)$
from the set of feasible solutions $X$ is pareto-optimal with respect to the criteria~$f$,
that does not really help the DM in decision making process.

In order to reduce the compromise set $P_f(X)$ let us consider the following information about the DM's preferences:
\begin{itemize}
  \item (P1) the group of the resources costs $\{f_1, f_2\}$ is more important then the income $\{f_3\}$
  with parameters $w_1^{(1)}, w_2^{(1)}, w_3^{(1)}$;
  \item (P2) the income $\{f_3\}$ is more important then the group of the resources costs $\{f_1, f_2\}$
  with parameters $w_1^{(2)}, w_2^{(2)}, w_3^{(2)}$.
\end{itemize}

First ''quantum of information'' (P1) says that the DM is ready to decrease the income $f_3$ by $w_3^{(1)}$ amount for
decreasing the costs of the basic production assets $f_1$ by $w_1^{(1)}$ amount and
the labor costs $f_2$ by $w_2^{(1)}$ amount simultaneously. And vice a versa, according to
second ''quantum of information'' (P2) the DM is willing to compromise
by increasing the costs of the basic production assets $f_1$ by $w_1^{(2)}$ amount and
the labor costs $f_2$ by $w_2^{(2)}$ amount simultaneously for gaining the additional income $w_3^{(2)}$.

Following the definition of ''quantum of information'' the DM's preferences (P1) and (P2) give the vectors
\begin{equation*}
y^{(1)} = ( w_1^{(1)}, w_2^{(1)}, -w_3^{(1)} ),
\quad
y^{(2)} = ( -w_1^{(2)}, -w_2^{(2)}, w_3^{(2)})
\end{equation*}
such that the relations $y^{(1)} \succ 0_3$ and $y^{(2)} \succ 0_3$ hold.

Obviously, statements (P1) and (P2) contradict each other. For that reason any collection of ''quantum of information''
should satisfy so-called {\it condition of consistency} (for the definition see~\cite{Nogin2016}).
With regards to ''quanta of information'' (P1) and (P2), it consists in implementation of at least one of the following inequalities:
\begin{equation}
\label{ref_W}
w_1^{(1)} / w_3^{(1)} > w_1^{(2)} / w_3^{(2)},
\quad
w_2^{(1)} / w_3^{(1)} > w_2^{(2)} / w_3^{(2)}.
\end{equation}
From natural point of view the compromise is justified if the gain is greater, than the loose.
Hence, item (P1) gives the inequalities $(w_1^{(1)} > w_3^{(1)})$ and $(w_2^{(1)} > w_3^{(1)})$, and
item (P2) establishes that $(w_3^{(2)} > w_1^{(2)})$ and $(w_3^{(2)} > w_2^{(2)})$. It leads to
implementation of both inequalities in~(\ref{ref_W}) simultaneously, that we will further assume.

\section{Pareto set reduction}
\label{ref_reduction}
Now apply information (P1) and (P2) to the mulictiteria choice model $<X, f, \succ>$.
From~\cite{Klimova_Noghin} we have the following theorem.

\begin{theorem}
Assume that there is a given information (P1) and (P2), and  both inequalities in~(\ref{ref_W}) hold simultaneously.
Then for any set of selected solutions $C(X)$
the inclusions $C(X)\subseteq P_g(X)\subseteq P_f(X)$ are valid.
Here, $P_g(X)$ is the Pareto set with respect to 4-dimensional vector criterion $g$ with components
\begin{equation*}
\begin{split}
\label{ref_g}
g_{13} = w_1^{(1)} f_3 + w_3^{(1)} f_1, & \quad
g_{23} = w_2^{(1)} f_3 + w_3^{(1)} f_2, \\
g_{31} = w_1^{(2)} f_3 + w_3^{(2)} f_1, & \quad
g_{32} = w_2^{(2)} f_3 + w_3^{(2)} f_2. \\
\end{split}
\end{equation*}
\end{theorem}

Thus, the use of information (P1) and (P2) consists in constructing the Pareto set of
a ''new'' multicriteria choice problem $<X, g, \succ>$. 
Taking into account formulae of functions $f_1$, $f_2$, and $f_3$~(\ref{ref_f}) we get the components of ''new'' vector criterion as functions of the basic production assets $K$ and the labor resources~$L$:
\begin{equation}
\begin{split}
\label{ref_g}
g_{13}(K,L) =& w_1^{(1)} p_{_Q} F(a K^{-r} + (1-a)L^{-r})^{\frac{1}{-r}} - w_3^{(1)} p_{_K} K, \\
g_{23}(K,L) =& w_2^{(1)} p_{_Q} F(a K^{-r} + (1-a)L^{-r})^{\frac{1}{-r}} - w_3^{(1)} p_{_L} L, \\
g_{31}(K,L) =& w_1^{(2)} p_{_Q} F(a K^{-r} + (1-a)L^{-r})^{\frac{1}{-r}} - w_3^{(2)} p_{_K} K, \\
g_{32}(K,L) =& w_2^{(2)} p_{_Q} F(a K^{-r} + (1-a)L^{-r})^{\frac{1}{-r}} - w_3^{(2)} p_{_L} L. \\
\end{split}
\end{equation}

The set of feasible solutions $X \equiv \R_+^2$ is convex.
Due to conditions~(\ref{ref_func_CES_cond}) functions $g_{13}(K,L)$, $g_{23}(K,L)$, $g_{31}(K,L)$, and
$g_{32}(K,L)$~(\ref{ref_g}) are concave on set $X$.
Thus, we obtain the problem with the convex set $X$ and the concave vector criterion $g$.
In this case, according to theorem of Karlin~\cite{Kaprlin} and Hurwicz~\cite{Arrow} the set of proper efficient points (which is slightly narrower than the Pareto set $P_g(X)$)
can be found by maximizing the linear combination of criteria components
\begin{multline}
\label{ref_phi}
\varphi(K, L) = \lambda_1 g_{_{13}}(K,L) + \lambda_2 g_{_{23}}(K,L) + \\
+ \lambda_3 g_{_{31}}(K,L) + \lambda_4 g_{_{32}}(K,L).
\end{multline}
with the positive coefficients $\lambda_1$, $\lambda_2$, $\lambda_3$, $\lambda_4$ such that $\Sigma_{i=1}^{4}\lambda_i = 1$.
If we substitute expressions~(\ref{ref_g}) in linear combination~(\ref{ref_phi}),
we get the following form of function $\varphi$:
\begin{equation*}
\varphi(K, L) = \bar{\lambda}_1 K + \bar{\lambda}_2 L + \bar{\lambda}_3 (a K^{-r} + (1-a)L^{-r})^{\frac{1}{-r}},
\end{equation*}
where quantities
\setlength{\arraycolsep}{0.0em}
\begin{eqnarray*}
\bar{\lambda}_1 = -p_{_K}(\lambda_{13}w_3^{(1)}&{}+{}&\lambda_{31}w_3^{(2)}), \nonumber\\
\bar{\lambda}_2 = -p_{_L}(\lambda_{23}w_3^{(1)}&{}+{}&\lambda_{32}w_3^{(2)}), \nonumber\\
\bar{\lambda}_3 = Fp_{_Q}(\lambda_{13}w_1^{(1)} + \lambda_{23}w_2^{(1)}&{}+{}&\lambda_{31}w_1^{(2)} + \lambda_{32}w_2^{(2)}). \nonumber
\end{eqnarray*}
\setlength{\arraycolsep}{5pt}

To find the extremum $x_0 = (K_0, L_0)$ of the function $\varphi(K, L)$ according to its necessity condition
evaluate first partial derivatives
\begin{equation*}
\varphi'_K = \bar{\lambda}_1 + \bar{\lambda}_3 a \left( a + (1-a)\left(\frac{L}{K}\right)^{-r}\right)^{\frac{1+r}{-r}},
\end{equation*}
\begin{equation*}
\varphi'_L = \bar{\lambda}_2 + \bar{\lambda}_3 (1-a) \left( a\left(\frac{K}{L}\right)^{-r} + (1-a)\right)^{\frac{1+r}{-r}},
\end{equation*}
and then solve the equations
\begin{equation}
\label{ref_varphi_eq}
\varphi'_K(K_0, L_0)~=~0, \ \ \varphi'_L(K_0, L_0)~=~0.
\end{equation}
One can check that the solution of equations~(\ref{ref_varphi_eq}) is a family of lines
\begin{equation}
\label{ref_sol_L0_K0}
L_0 = K_0 \left(\left(\frac{(1-a)\bar{\lambda}_1}{a \bar{\lambda}_2}\right)^{\frac{-r}{1+r}}-
\frac{a}{1-a} \right)^{\frac{-1}{r}}.
\end{equation}

Due to the concavity of functions $g_{13}(K,L)$, $g_{23}(K,L)$, $g_{31}(K,L)$, and $g_{32}(K,L)$~(\ref{ref_g}) on set $X$
the linear combination $\varphi(K, L)$ also will be concave on set $X$.
This means that for any $(K, L) \in X$ the Hessian
$|H(\varphi(K, L))|=$ $=\varphi''_{_{K^2}}(K, L)\varphi''_{_{L^2}}(K, L) - (\varphi''_{_{KL}}(K, L))^2 < 0$ including the point $(K_0, L_0)$.
So, the sufficient condition of maximum is justified, and we get that point $(K_0, L_0)$ is a maximum point of
the function $\varphi(K, L)$.

Thus, we have the following set of compromise (the Pareto set) of the muticriteria choice problem with criteria $g$:
\begin{equation}
\label{ref_P_g}
P_g(X) = \{(K_0,L_0) \in X \mid (\ref{ref_sol_L0_K0}) \ \mbox{holds} \}.
\end{equation}

\section{Fuzzy multicriteria choice problem}
\subsection{Mathematical background}
Following previous sections the DM could arrive to one of three possibilities when considering
two feasible solutions $x', \ x'' \in X$:
the first $x'$ solution is preferable, the second solution $x''$ is preferable, and two solutions are non-comparable.
At the same time it seems more natural to consider such preference relation, which also reflects a wishes degree of confidence.
And we come to the notation of fuzzy preference relation generalizing aforementioned preference relation $\succ$.
In~\cite{Nogin2016} author states a fuzzy multicriteria choice problem, including corresponding axioms of ''reasonable''
choice and the {\it Edgeworth--Pareto principle} (both in terms of fuzzy case), which will be given below.

First recall the definitions of fuzzy set and fuzzy binary relation~\cite{Zadeh}.
Let $U$ be some set, that we call universe.

\begin{definition}
A fuzzy set $A$ in $U$ is characterized by a membership function $f_A(x)$ which associates each point in $U$
a real number in interval [0, 1], with the value of $f_A(x)$ at $x$ representing the ''grade of membership'' of $x$ in $A$.
\end{definition}

\begin{definition}
A fuzzy relation in $U$ is a fuzzy set in the product space $U \times U$.
\end{definition}

Upon the latter let us introduce a fuzzy preference relation of the DM
with a membership function $\mu \colon Y \times Y \to [0, 1]$ as follows.
If for any solutions $x', x'' \in X$ the equality $\mu (f(x'), f(x''))=\mu^\ast$ holds,
then the DM prefers the solution $x'$ to the solution $x''$ with degree of confidence $\mu^\ast$
showing how much {\color{black}it is sure} in the choice. 

Thus, we get the following fuzzy multicriteria choice problem $<X, f, \mu>$~\cite{Nogin2016}:
\begin{itemize}
\item a set of feasible solutions (alternatives) $X$;
\item a vector criterion $f = (f_1, f_2, \ldots,  f_m)$ defined on set $X$;
\item a fuzzy preference relation of the DM $\mu$ defined on set $Y = f(X)$.
\end{itemize}

Analogously subsection~\ref{ref_subsect_crisp_model_math} we suppose, that the fuzzy relation $\mu$ satisfies
the axioms of ''fuzzy reasonable'' choice generalizing corresponding axioms of crisp case~\cite{Nogin2016}.
So, the fuzzy relation $\mu$ is irreflexive, transitive,
invariant under linear positive transformation and
compatible with each criteria $f_1, f_2, \ldots, f_m$ in terms of fuzzy sets.
The {\it Edgeworth--Pareto principle} has the following form: under class of axioms of ''fuzzy reasonable'' choice
the inequality $\lambda_X^C(x) \leqslant \lambda_X^P(x)$ is valid for all $x \in X$.
Here, $\lambda_X^P(\cdot)$ is a membership function of the Pareto set (crisp) $P_f(X)$,
and $\lambda_X^C(\cdot)$ is a membership function of the set of selected solutions $C(X)$.

This principle specifies an upper bound of ''fuzzy'' choice and restricts the class of fuzzy multicriteria choice problem,
in which the choice should be done among the Pareto set, and only among it.
Under ''fuzzy'' choice, fuzzy set of selected solutions we imply some hypothetical fuzzy set,
that meets all possible preferences of the DM.

Pass to the generalization of Definition~\ref{ref_def_AB}.
\begin{definition}
\label{ref_AB_fuzzy}
Say that there exists a {\it ''fuzzy quantum of information''} about the DM's preference
relation $\mu$, if the vector $y' \in \R^m$ with components~(\ref{ref_AB})
satisfies the expression $\mu(y', 0_m) = \mu^\ast$, where value $\mu^\ast$ belongs to interval $[0, 1]$
and shows the degree of confidence of compromise: the group of criteria $A$ is more important than the group of criteria $B$ with given positive parameters $w_i \ \forall i \in A$, $w_j \ \forall j \in B$.
\end{definition}

Using some collection of ''fuzzy quantum of information'' we could derive an upper bound
of the ''fuzzy'' optimal choice (fuzzy set of selected vectors).
According to~\cite{Nogin2016} it consists in solving series of crisp multicriteria choice problem.
We consider such process more precisely in the next subsection.

\subsection{Fuzzy economic problem}
\label{ref_fuzzy_model}
\subsubsection{Problem statement}
Return to the multicriteria economic problem considered in subsection~\ref{ref_econ_sys_crysp}
with the set of feasible resources $X \equiv \R^2_{+}$, the vector criterion $f = (f_1, f_2, f_3)$ defined by~(\ref{ref_f}),
and the DM's preference relation $\succ$ on set $Y = f(X)$.

Now pass to the ''fuzzy'' problem replacing crisp DM's preference relation $\succ$ by the fuzzy relation $\mu$.
All axioms of ''fuzzy reasonable'' choice and the Edgeworth-Pareto principle are valid for this problem. 
So, the Pareto set $P_f(X) = X$ is an upper bound of optimal choice, and we cannot make it tighter unless
taking into account some additional information in terms of Definition~\ref{ref_AB_fuzzy}, like in the crisp case.

Following Definition~\ref{ref_AB_fuzzy} let us extend information (P1) and (P2) to the relation $\mu$
and indicate it by (FP1) and (FP2).
\begin{itemize}
  \item (FP1) the group of the resources costs $\{f_1, f_2\}$ is more important then the income $\{f_3\}$
  with parameters $w_1^{(1)}, w_2^{(1)}, w_3^{(1)}$ and degree of confidence $\mu_1$;
  \item (FP2) the income $\{f_3\}$ is more important then the group of the resources costs $\{f_1, f_2\}$
  with parameters $w_1^{(2)}, w_2^{(2)}, w_3^{(2)}$ and degree of confidence $\mu_2$.
\end{itemize}
Next consider how to use such information to reduce the compromise set $P_f(X)$.
Note, that condition~(\ref{ref_W}) is sufficient for the consistency of information (FP1) and (FP2).
Moreover, we suppose that both inequalities in~(\ref{ref_W}) are valid due to natural requirements
(see details in subsection~\ref{ref_econ_sys_crysp}).

\subsubsection{Applying the information}
Let the inequality $\mu_1 \geqslant \mu_2$ is true. According to the results in~\cite{Klimova} the use of
information (FP1) and (FP2) consists in solving three crisp multicriteria problems,
that means finding a corresponding Pareto set.
This process yields the fuzzy set, which is an upper bound of the ''fuzzy'' optimal choice.
By $\lambda_X^M(\cdot)$ we further denote the membership function of this set.

Firstly we should solve the problem $<X, f, \succ>$ when we have not any additional information
in terms of Definition~\ref{ref_AB}, i.e. we should find the Pareto set $P_f(X)$.
Actually, we have already done it, $P_f(X) = X$. And then put $\lambda_X^M(x) = 1$
for all solutions $x \in P_f(X)$, and $\lambda_X^M(x) = 0$ for all solutions $x \in X \setminus P_f(X)$.

Secondly we should solve the problem $<X, \bar{f}, \succ>$ with information (P1),
where the vector criterion $\bar{f} = (\bar{f}_1, \bar{f}_2, \bar{f}_3, \bar{f}_4)$
has the following components: $\bar{f}_1 = f_1$, $\bar{f}_2 = f_2$, $\bar{f}_3 = g_{13}$, and $\bar{f}_4 = g_{23}$.
Here, the functions $f_1$ and $f_2$ is defined by~(\ref{ref_f}),
and the functions $g_{13}$ and $g_{23}$ is defined by~(\ref{ref_g}).
And then put $\lambda_X^M(x) = 1 - \mu_1$ for all solutions $x \in P_f(X) \setminus P_{\bar{f}}(X)$.

Thirdly we will find a solution to the problem $<X, g, \succ>$ with information (P1) and (P2), where the vector criterion $g$
defined according to formulae~(\ref{ref_g}). Next we should put $\lambda_X^M(x) = 1 - \mu_2$
for all solutions $x \in P_{\bar{f}}(X) \setminus P_g(X)$.
And any solution $x$ from set $P_g(X)$ still has degree of confidence equaled 1.

In~\cite{Klimova} author proved that the following inclusions are valid:
\begin{equation*}
C(X) \subseteq P_g(X) \subseteq P_{\bar{f}}(X) \subseteq P_f(X),
\end{equation*}
or in terms of membership functions we have
\begin{equation}
\label{ref_lambda_M}
\lambda_X^C(x) \leqslant \lambda_X^M(x) \leqslant \lambda_X^P(x) \quad \forall x \in X.
\end{equation}
So, aforementioned subtractions of corresponding sets are not unreasonable.

Thus, solutions of three ''crisp'' multicriteria problems give a fuzzy set with membership function $\lambda_X^M(\cdot)$
such that~(\ref{ref_lambda_M}). Now, consider them more precisely.

\subsubsection{Multicriteria problem $<X, \bar{f}, \succ>$}
We have the set of feasible solutions $X$ and the vector criterion $\bar{f}~=~(\bar{f}_1, \bar{f}_2, \bar{f}_3, \bar{f}_4)$.
Obviously, the functions $\bar{f}_1, \bar{f}_2, \bar{f}_3, \bar{f}_4$ are concave.
And thus, as mentioned in section~\ref{ref_reduction},
the set of proper efficient points (which is slightly narrower than the Pareto set)
can be found by maximizing the linear combination of criteria components
\begin{multline*}
\label{ref_phi_fuzzy1}
\bar{\varphi}(K, L) = \lambda_1 f_1(K,L) + \lambda_2 f_2(K,L) + \\
+ \lambda_3 g_{_{13}}(K,L) + \lambda_4 g_{_{23}}(K,L).
\end{multline*}
with the positive coefficients $\lambda_1$, $\lambda_2$, $\lambda_3$, $\lambda_4$ such that $\Sigma_{i=1}^{4}\lambda_i = 1$.
Substituting the corresponding expressions we get
\setlength{\arraycolsep}{0.0em}
\begin{eqnarray*}
\bar{\varphi}(K, L) = - p_{_K}(\lambda_1 + \lambda_3 w_3^{(1)})K - p_{_L}(\lambda_2 + \lambda_4 w_3^{(1)})L\:{+}\nonumber\\
{+}\:p_{_Q}F(\lambda_3 w_1^{(1)} + \lambda_4 w_2^{(1)}) (a K^{-r} + (1-a)L^{-r})^{-(1/r)}.\nonumber\\
\end{eqnarray*}
\setlength{\arraycolsep}{5pt}
According to necessity condition of maximum evaluate first partial derivatives of the function $\bar{\varphi}(K, L)$:
\setlength{\arraycolsep}{0.0em}
\begin{eqnarray*}
\bar{\varphi}'_K&{}={}&- p_{_K}(\lambda_2 + \lambda_4 w_3^{(1)}) + a p_{_Q}F(\lambda_3 w_1^{(1)} + \lambda_4 w_2^{(1)})\:{\times}\nonumber\\
&&\:{\times}\left( a + (1-a)\left(\frac{L}{K}\right)^{-r}\right)^{\frac{1+r}{-r}},\nonumber\\
\bar{\varphi}'_L&{}={}&- p_{_L}(\lambda_1 + \lambda_3 w_3^{(1)}) + (1-a) p_{_Q}F(\lambda_3 w_1^{(1)} + \lambda_4 w_2^{(1)})\:{\times}\nonumber\\
&&\:{\times}\left( a\left(\frac{K}{L}\right)^{-r} + (1-a)\right)^{\frac{1+r}{-r}},\nonumber
\end{eqnarray*}
\setlength{\arraycolsep}{5pt}
And then solving the equations $\bar{\varphi}'_K(K_0, L_0)~=~0$, $\bar{\varphi}'_L(K_0, L_0)~=~0$
we have a linear dependence
\begin{equation}
\label{ref_bar_f_L0_K0}
L_0 = K_0 \left(\left(\frac{(1-a)p_{_K}(\lambda_2 + \lambda_4 w_3^{(1)})}
{a p_{_L}(\lambda_1 + \lambda_3 w_3^{(1)})}\right)^{\frac{-r}{1+r}} - \frac{a}{1-a} \right)^{\frac{-1}{r}}.
\end{equation}

Due to the concavity of functions $f_1(K,L)$, $f_2(K,L)$, $g_{13}(K,L)$, and $g_{23}(K,L)$ on set $X$
the linear combination $\bar{\varphi}(K, L)$ also will be concave on set $X$.
So, the sufficient condition of maximum is justified: the Hessian $|H(\bar{\varphi}(K_0, L_0))|~<~0$.

Thus, we have the Pareto set
\begin{equation}
\label{ref_P_bar_f}
P_{\bar{f}}(X) = \{(K_0,L_0) \in X \mid (\ref{ref_bar_f_L0_K0}) \mbox{ holds} \}
\end{equation}
of the multicriteria choice problem $<X, \bar{f}, \succ>$. Then values of the membership function $\lambda_X^M(\cdot)$
for all solutions $x~\in~P_f(X) \setminus P_{\bar{f}}(X)$ is $1 - \mu_1$.

\subsubsection{Constructing membership function $\lambda_X^M(\cdot)$}
\label{sect_lambda_M}
As we mentioned before in this subsection the first multicriteria problem gives the Pareto set $P_f(X) = X$,
and we have $\lambda_X^M(x) = 1$ $\forall x \in X$.

From the second multicriteria problem $<X, \bar{f}, \succ>$ we get $\lambda_X^M(x) = 1 - \mu_1$
$\forall x \in X \setminus P_{\bar{f}}(X)$, where $P_{\bar{f}}(X)$ defined by~(\ref{ref_P_bar_f}).

Finally, the third multicriteria problem $<X, g, \succ>$ has been actually solved in section~\ref{ref_reduction},
and we establish $\lambda_X^M(x) = 1 - \mu_2$ $\forall x \in P_{\bar{f}}(X) \setminus P_g(X)$,
where $P_g(X)$ defined by~(\ref{ref_P_g}). Note, that for solutions $x \in P_g(X)$ the equality
$\lambda_X^M(x) = 1$ holds.

Thus, we construct the fuzzy set with membership function $\lambda_X^M(\cdot)$, which is an upper bound of
the fuzzy set of selected solutions $C(X)$, and this bound belongs to the initial Pareto set $P_f(X)$.
Let us estimate the subtractions of sets $P_f(X) \setminus P_{\bar{f}}(X)$ and $P_{\bar{f}}(X) \setminus P_g(X)$.

Start with defining the condition when a pair $x'~=~(K',L')~\in~P_f(X)$. Since the components of
vector criterion $f$ are concave,
a pair $x'$ is pareto-optimal if and only if it is a maximum point of function
\setlength{\arraycolsep}{0.0em}
\begin{eqnarray*}
\label{ref_phi_0}
&\varphi_0(K, L)=\lambda_{01} f_1(K,L) + \lambda_{02} f_2(K,L) + \lambda_{03} f_3(K,L)\:= \\
&\:=-\lambda_{01} p_{_K} K - \lambda_{02} p_{_L} L + \lambda_{03} p_{_Q} F(a K^{-r} + (1-a)L^{-r})^\frac{-1}{r}.
\end{eqnarray*}
\setlength{\arraycolsep}{5pt}
where the coefficients $\lambda_{01}$, $\lambda_{02}$, and $\lambda_{03}$ are positive and
such that $\Sigma_{i=1}^{3}\lambda_{0i} = 1$.
Then, analogously the previous statements we get that point $x'$ satisfies the following linear dependence
\begin{equation}
\label{ref_eq_P_bar_f}
L' = K' \left(\left(\frac{(1-a)p_{_K}\lambda_{01}}{a p_{_L}\lambda_{02}}\right)^{\frac{-r}{1+r}} -
\frac{a}{1-a} \right)^{\frac{-1}{r}}.
\end{equation}
At the same time, the condition $x' \notin P_{\bar{f}}(X)$ should be valid,
which means that expression~(\ref{ref_bar_f_L0_K0}) is not true for $K_0 = K'$, $L_0 = L'$.
Thus, the inclusion $x' \in P_f(X) \setminus P_{\bar{f}}(X)$ is equivalent
to equation~(\ref{ref_eq_P_bar_f}), and there is no such positive quantities
$\lambda_1$, $\lambda_2$, $\lambda_3$, and $\lambda_4$ that
$(\lambda_1 + w_3^{(1)}\lambda_2)\lambda_{01} =$ $= (\lambda_2 + \lambda_4 w_3^{(1)})\lambda_{02}$.

Following similar arguments we can come to the estimation of $P_{\bar{f}}(X) \setminus P_g(X)$.
The inclusion $x_0 \in P_{\bar{f}}(X) \setminus P_g(X)$ is equivalent
to equation~(\ref{ref_bar_f_L0_K0}), and there is no such positive quantities
$\lambda_{13}$, $\lambda_{23}$, $\lambda_{31}$, and $\lambda_{32}$ that
\begin{multline*}
(\lambda_1 + \lambda_3 w_3^{(2)})(\lambda_{13}w_3^{(1)} + \lambda_{31}w_3^{(2)}) = \\
= (\lambda_2 + \lambda_4 w_3^{(1)})(\lambda_{23}w_3^{(1)} + \lambda_{32}w_3^{(2)}).
\end{multline*}

\subsubsection{Case $\mu_1 < \mu_2$}
Now consider the case, when the inequality $\mu_1 < \mu_2$ is valid. Then we have the following changes
in the results derived above.

The second multicriteria problem is $<X, \hat{f}, \succ>$ instead of $<X, \bar{f}, \succ>$,
where the vector criterion $\hat{f} = (\hat{f}_1, \hat{f}_2, \hat{f}_3)$
has the following components: $\hat{f}_1 = g_{31}$, $\hat{f}_2 = g_{32}$, $\hat{f}_3 = f_3$.
Here, the function $f_3$ is defined by~(\ref{ref_f}), and the functions $g_{31}$ and $g_{32}$ are defined by~(\ref{ref_g}).
Hence, for all solutions $x \in P_f(X) \setminus P_{\hat{f}}(X)$
we should put $\lambda_X^M(x) = 1 - \mu_2$.
And the expression $\lambda_X^M(x) = 1 - \mu_1$ is valid
for all solutions $x \in P_{\hat{f}}(X) \setminus P_g(X)$.





\section*{Acknowledgment}
This research is supported by the Russian Foundation for Basic Research grant 17-07-00371.




%

\end{document}